\documentclass[11pt,letterpaper,reqno]{amsart}
\usepackage{amsmath,amsthm,amsfonts,amssymb,mathtools,algpseudocode}
\usepackage{comment,bookmark,bm,color}

\hypersetup{pdfstartview={FitH}}

\newtheorem{Theorem}{{\bf Theorem}}[section]

\newtheorem{Algorithm}[Theorem]{{\bf Algorithm}}
\newtheorem{Proposition}[Theorem]{{\bf Proposition}}

\numberwithin{equation}{section}

\newcommand{\calA}{\mathcal{A}}
\newcommand{\calB}{\mathcal{B}}
\newcommand{\calC}{\mathcal{C}}
\newcommand{\calE}{\mathcal{E}}
\newcommand{\calG}{\mathcal{G}}

\newcommand{\calT}{\mathcal{T}}

\newcommand{\R}{\mathbb{R}}
\newcommand{\N}{\mathbb{N}}
\newcommand{\C}{\mathbb{C}}
\newcommand{\anti}{\text{anti}}
\newcommand{\sign}{\text{sign}}
\newcommand{\ran}{\text{rank}}
\newcommand{\de}{\text{det}}

\begin{document}

\title[CP decomposition and low-rank approximation of antisymmetric tensors]{CP decomposition and low-rank approximation of antisymmetric tensors}
%\author{One author\thanks{some info} \and Another author\thanks{more info}}
\author{Erna Begovi\'{c}~Kova\v{c}}\thanks{\textsc{Erna Begovi\'{c} Kova\v{c}}, University of Zagreb Faculty of Chemical Engineering and Technology, Maruli\'{c}ev trg 19, 10000 Zagreb, Croatia, \texttt{ebegovic@fkit.hr}}
\author{Lana Peri\v{s}a}\thanks{\textsc{Lana Peri\v{s}a}, Visage Technologies, Ivana Lučića 2a, 10000 Zagreb, Croatia, \texttt{lana.perisa@visagetechnologies.com}}

\thanks{This work has been supported in part by Croatian Science Foundation under the project UIP-2019-04-5200.}
\date{\today}

\renewcommand{\subjclassname}{\textup{2020} Mathematics Subject Classification}

\begin{abstract}
For the antisymmetric tensors the paper examines a low-rank approximation which is represented via only three vectors.
We describe a suitable low-rank format and propose an alternating least squares structure-preserving algorithm for finding such approximation.
Moreover, we show that this approximation problem is equivalent to the problem of finding the best multilinear low-rank antisymmetric approximation and, consequently, equivalent to the problem of finding the best unstructured rank-$1$ approximation.
The case of partial antisymmetry is also discussed.
The algorithms are implemented in Julia programming language and their numerical performance is discussed.
\end{abstract}

\maketitle

\section{Introduction}

Tensor decompositions have been extensively studied in recent years~\cite{AcarSurvey,DeLMV_HOSVD,DeLMV04,HKD20,KB09}. However, the research was mostly focused on either unstructured or symmetric~\cite{CGLMsym08,Kolda15} tensors. In this paper we explore the antisymmetric tensors, their CP decomposition and the algorithms for the low-rank approximation.

The idea of the CP decomposition is to write a tensor as a sum of its rank-one components. It was first introduced by Hitchcock~\cite{Hitch27a,Hitch27b} in 1927, but it only became popular in the 1970s as CANDECOMP (canonical decomposition)~\cite{candecomp70} and PARAFAC (parallel factors)~\cite{parafac70}. This decomposition is closely related to the tensor rank $R$, which is defined as the minimal number of rank-$1$ summands in the exact CP decomposition. Contrary to the matrix case, the rank of a tensor can exceed its dimension, and it can be different over $\R$ and over $\C$. It is known that the problem of finding the rank of a given tensor is NP-hard.

When computing the CP approximation, the main question is the choice of the number of rank-one components. Given the antisymmetric structure of our tensors in question, we impose an additional constraint on the CP decomposition. This constraint assures that the resulting tensor is, indeed, antisymmetric, and it gives a bound on the minimal number of rank-one components.

We focus on the tensors of order three. For a given antisymmetric tensor $\calA\in\R^{n\times n\times n}$ our goal is to find its low-rank antisymmetric approximation which is represented via only three vectors. In particular, we are looking for the approximation $\widetilde{\calA}$ of $\calA$ such that $\ran(\widetilde{\calA})\leq6$ for any $n$, and
$$\widetilde{\calA}=\frac{1}{6}(x\circ y\circ z+y\circ z\circ x+z\circ x\circ y-x\circ z\circ y-y\circ x\circ z-z\circ y\circ x),$$
where $x,y,z\in\R^n$.
We propose an alternating least squares structure-preserving algorithm for solving this problem. The algorithm is based on solving a minimization problem in each tensor mode. We compare our algorithm with a ``naive'' idea which uses a posteriori antisymmetrization. Further on, we show that our approximation problem is equivalent to the problem of the best multilinear rank-$3$ structure-preserving antisymmetric tensor approximation from~\cite{BeKre17} and, consequently, equivalent to the problem of the best unstructured rank-$1$ approximation. This establishes the equivalence between our algorithm and the higher-order power method (HOPM). Therefore, corresponding convergence result for HOPM from~\cite{Usch15} can be applied.

Additionally, we study the tensors with partial antisymmetry, that is, antisymmetry in only two modes. Similarly to what we do for the tensors that are antisymmetric in all modes, we first determine a suited format of the CP decomposition, which is going to be simpler for the partial antisymmetry. Based on this format, for a given tensor $\calC\in\R^{n\times n\times m}$ antisymmetric in two modes, we are looking for its approximation $\widetilde{\calC}$ of the same structure such that $\widetilde{\calC}$ is represented by three vectors and $\ran(\widetilde{\calC})=2$.

In Section~\ref{sec:pre} we introduce the notation and preliminaries. Our problem of antisymmetric tensor approximation is described in Section~\ref{sec:problem}. In Sections~\ref{sec:algorithm1} we describe the approach with a posteriori antisymmetrization, while in Section~\ref{sec:algorithm2} we propose the algorithm for solving the minimization problem from Section~\ref{sec:problem}. Section~\ref{sec:partial} deals with the case of partial antisymmetry. In Section~\ref{sec:num} we discuss our numerical results obtained in Julia programming language and conclusion is given in Section~\ref{sec:conclusion}.

\section{Notation and preliminaries}\label{sec:pre}

Throughout the paper we denote tensors by calligraphic letters, e.g., $\calA$. We refer to the tensor dimension as its \emph{order}. Then, for $\calA\in\R^{n_1\times n_2\times\cdots\times n_d}$ we say that $\calA$ is a tensor of order $d$. Tensor $\calA\in\R^{n_1\times n_2\times\cdots\times n_d}$ is \emph{cubical} if $n_1=n_2=\cdots=n_d$.
Vectors obtained from a tensor by fixing all indices but the $m$th one are called \emph{mode-$m$ fibers}. Fibers of an order-$3$ tensor are columns (mode-$1$ fibers), rows (mode-$2$ fibers), and tubes (mode-$3$ fibers). Matrices obtained from a tensor by fixing all indices but two are called \emph{slices}.
Matrix representation of a tensor $\calA\in\R^{n_1\times n_2\times\cdots\times n_d}$ is called mode-$m$ \emph{matricization} and it is denoted by $A_{(m)}$. It is obtained by arranging mode-$m$ fibers of $\calA$ as columns of $A_{(m)}$.
The \emph{mode-$m$ product} of a tensor $\calA\in\R^{n_1\times n_2\times\cdots\times n_d}$ with a matrix $M\in\R^{p\times n_m}$ is a tensor
$\calB\in\R^{n_1\times\cdots\times n_{m-1}\times p\times n_{m+1}\times\cdots\times n_d}$,
$$\calB=\calA\times_m M, \quad \text{such that} \quad B_{(m)}=M A_{(m)}.$$
Tensor \emph{norm} is a generalization of the Frobeius norm. For $\calA\in\R^{n_1\times n_2\times\cdots\times n_d}$ we have
$$\|\calA\|=\sqrt{\sum_{i_1=1}^{n_1}\sum_{i_2=1}^{n_2}\cdots\sum_{i_d=1}^{n_d} a_{i_1i_2\ldots i_d}^2}.$$
The \emph{inner product} of two tensors $\calA,\calB\in\R^{n_1\times n_2\times\cdots\times n_d}$ is given by
$$\langle\calA,\calB\rangle=\sum_{i_1=1}^{n_1}\sum_{i_2=1}^{n_2}\cdots\sum_{i_d=1}^{n_d} a_{i_1i_2\ldots i_d}b_{i_1i_2\ldots i_d}.$$
The vector \emph{outer product} is denoted by $\circ$. A tensor $\calA\in\R^{n_1\times n_2\times\cdots\times n_d}$ is a rank-$1$ tensor if it can be written as the outer product of $d$ vectors,
$$\calA=v^{(1)}\circ v^{(2)}\circ\cdots\circ v^{(d)}.$$
Then
$$a_{i_1i_2\cdots i_d}=v^{(1)}_{i_1}v^{(2)}_{i_2}\cdots v^{(d)}_{i_d}, \quad 1\leq i_k\leq n_k, \ 1\leq k\leq d.$$

The \emph{Khatri-Rao product} of two matrices $A\in\R^{m\times n}, B\in \R^{p\times n}$ is defined as
$$A\odot B = \begin{bmatrix}
         a_1 \otimes b_1 & a_2 \otimes b_2 & \cdots & a_n \otimes b_n
        \end{bmatrix} \in\R^{(mp)\times n},$$
where $a_k$ and $b_k$ denote the $k$th column of $A$ and $B$, respectively.
The \emph{Hadamard (element-wise) product} of two matrices $A,B\in\R^{m\times n}$, is defined as
$$A\ast B=\begin{bmatrix}
          a_{11}b_{11} & a_{12}b_{12} & \cdots &a_{1n}b_{1n} \\
          a_{21}b_{21} & a_{22}b_{22} & \cdots &a_{2n}b_{2n} \\
          \vdots & \vdots & \ddots & \vdots\\
          a_{m1}b_{m1} & a_{m2}b_{m2} & \cdots &a_{mn}b_{mn}
         \end{bmatrix}\in\R^{m\times n}.$$
The \emph{Moore-Penrose inverse} of $A$ is denoted by $A^+$.

For a tensor $\calA\in\R^{n_1\times n_2\times n_3}$, its \emph{CP approximation} takes the form
\begin{equation}\label{cpvec}
\calA\approx\sum_{i=1}^r (x_i\circ y_i\circ z_i),
\end{equation}
where $x_i\in\R^{n_1}$, $y_i\in\R^{n_2}$, $z_i\in\R^{n_3}$.
If we arrange vectors $x_i,y_i,z_i$, $i=1,\ldots,r$ into matrices
$$X=\left[
      \begin{array}{cccc}
        x_1 & x_2 & \cdots & x_r \\
      \end{array}
    \right], \quad
Y=\left[
      \begin{array}{cccc}
        y_1 & y_2 & \cdots & y_r \\
      \end{array}
    \right], \quad
Z=\left[
      \begin{array}{cccc}
        z_1 & z_2 & \cdots & z_r \\
      \end{array}
    \right],
$$
relation~\eqref{cpvec} can be written as
\begin{equation}\label{cp}
\calA\approx[[X,Y,Z]]=\sum_{i=1}^r (x_i\circ y_i\circ z_i).
\end{equation}
The smallest number $r$ in the exact CP decomposition~\eqref{cp} is called \emph{tensor rank}. We write $\ran(\calA)=r$.

The most commonly used algorithm for computing the CP approximation is alternating least squares (ALS) algorithm. (See, e.g.,\@ \cite{KB09}.) In Algorithm~\ref{agm:cpals} we give the CP-ALS algorithm for order-$3$ tensors.
\bigskip

\begin{Algorithm}\label{agm:cpals}
\vspace{0.5ex}\hrule\vspace{0.5ex}
\emph{CP-ALS}
\vspace{0.5ex}\hrule
\begin{algorithmic}
\State \textbf{Input:} $\calA\in\R^{n\times n\times n}$, $r\in\N$
\State \textbf{Output:} $X,Y,Z\in\R^{n\times r}$
\State Initialize $X,Y,Z$ as leading $r$ left singular vectors of $A_{(i)}$, $i=1,2,3$, respectively.
\Repeat
\State $X=A_{(1)}(Z\odot Y)(Y^TY\ast Z^TZ)^+$
\State $Y=A_{(2)}(Z\odot X)(X^TX\ast Z^TZ)^+$
\State $Z=A_{(3)}(Y\odot X)(X^TX\ast Y^TY)^+$
\Until{convergence or maximum number of iterations}
\end{algorithmic}
\hrule
\end{Algorithm}
\bigskip

\section{Problem description}\label{sec:problem}

A cubical tensor is \emph{symmetric} (sometimes also called supersymmetric) if its elements are invariant to any permutation of indices.
On the contrary, a cubical tensor is \emph{antisymmetric} if its elements change sign when permuting pairs of indices. In particular, an order-$3$ tensor $\calA\in\R^{n\times n\times n}$ is antisymmetric if
\begin{equation}\label{antisymmetry}
a_{ijk}=a_{jki}=a_{kij}=-a_{ikj}=-a_{jik}=-a_{kji}, \quad 1\leq i,j,k\leq n.
\end{equation}
Such tensors are also called alternating $3$-tensors $\Lambda^3(\R_n)$~\cite{Lan12}, or $3$-vectors~\cite{Greub}.
The antisymmetric tensors appear in applications such as quantum chemistry~\cite{quantum} and electromagnetism~\cite{electromagnetism}. Besides, they are interesting from the mathematical point of view~\cite{BeKre17,Hackbusch18}.
From the definition of the antisymmetric tensor $\calA$ it obviously follows that
\begin{itemize}
\item[(i)] In all modes, all slices of $\calA$ are antisymmetric matrices.
\item[(ii)] In all modes, all slices have one null column and one null row.
\item[(iii)] Antisymmetric tensors are data-sparse in the sense that many of its non-zero elements are the same, up to the sign.
\end{itemize}
These facts are useful when it comes to the implementation of the specific algorithms.

We can define the \emph{antisymmetrizer} ``$\anti$'' as the orthogonal projection of a general cubical order-$d$ tensor $\calB$ to the subspace of antisymmetric tensors. Then, $\calA=\anti(\calB)$ is an order-$d$ tensor given by
$$\calA(i_1,i_2,\ldots,i_d)\coloneqq\frac{1}{d!}\sum_{p\in\pi(d)}\sign(p)\calB(p(i_1),p(i_2),\ldots,p(i_d)),$$
where $\pi(d)$ denotes the set of all permutations of length $d$. Hence, for $d=3$, $\calB\in\R^{n\times n\times n}$ and $\calA=\anti(\calB)$ we have
\begin{equation}\label{anti3d}
a_{ijk}=\frac{1}{6}(b_{ijk}+b_{jki}+b_{kij}-b_{ikj}-b_{jik}-b_{kji}).
\end{equation}

Let $\calA\in\R^{n\times n\times n}$ be an antisymmetric tensor of order three. Take a triplet of indices $(i,j,k)$, $1\leq i<j<k\leq n$. It follows from~\eqref{antisymmetry} that a subtensor $\hat{\calA}$ of $\calA$ obtained at the intersection of the $i$th, $j$th, and $k$th column, row, and tube is of the form
$$\hat{\calA}=\alpha\calE,$$
where $\alpha\in\R$ and $\calE$ is a $3\times3\times3$ tensor such that
$$\calE(i_1,i_2,i_3)=\left\{
                       \begin{array}{rl}
                         1, & \text{if the indices make an even permutation of $(1,2,3)$,} \\
                         -1, & \text{if the indices make an odd permutation of $(1,2,3)$,} \\
                         0, & \text{if two or more indices are equal.} \\
                       \end{array}
                     \right.
$$
Tensor $\calE$ is called the \emph{Levi-Civita tensor}~\cite{GoldsteinCM}.
We can also write $\calE$ using its matricization
\begin{equation}\label{tensorE}
E_{(1)}=\left[
            \begin{array}{ccc|ccc|ccc}
              0 & 0 & 0 & 0 & 0 & -1 & 0 & 1 & 0 \\
              0 & 0 & 1 & 0 & 0 & 0 & -1 & 0 & 0 \\
              0 & -1 & 0 & 1 & 0 & 0 & 0 & 0 & 0 \\
            \end{array}
          \right].
\end{equation}
Obviously, $\calE$ is the simplest possible antisymmetric non-zero order-$3$ tensor.

For three given vectors $x,y,z\in\R^n$ we define an $n\times n\times n$ antisymmetric tensor associated to these vectors as
\begin{equation}\label{antiform}
\calA_6(x,y,z)\coloneqq\frac{1}{6}(x\circ y\circ z+y\circ z\circ x+z\circ x\circ y-x\circ z\circ y-y\circ x\circ z-z\circ y\circ x).
\end{equation}
Note that tensor $\calE$ is a special case of the antisymmetric tensor $\calA_6(x,y,z)$.
For $x=[6,0,0]^T$, $y=[0,1,0]^T$, $z=[0,0,1]^T$, we get $\calA_6(x,y,z)=\calE$.
Moreover, for a rank one tensor $\calT=[[x,y,z]]$, we have $\calA_6(x,y,z)=\anti(\calT)$.
Tensor format~\eqref{antiform} can be favourable because it represents an antisymmetric tensor via only three vectors, that is, $3n$ entries.
On the other hand, the standard form of an $n\times n\times n$ antisymmetric tensor contains $\binom{n}{3}$ different entries.
Besides, tensor $\calA_6(x,y,z)$ is a low-rank tensor. For any size $n$, we have $\ran(\calA_6(x,y,z))\leq6$.

Our goal is to approximate a given antisymmetric tensor $\calA$ with a low-rank antisymmetric tensor of the form~\eqref{antiform}. We demonstrate two approaches. The ``naive'' one is given in Section~\ref{sec:algorithm1}. Then, in Section~\ref{sec:algorithm2} we formulate this problem as the minimization probelm. For a given non-zero antisymmetric tensor $\calA\in\mathbb{R}^{n\times n\times n}$, we are looking for a tensor $\widetilde{\calA}=\calA_6(x,y,z)$, i.e., vectors $x,y,z\in\R^n$, such that
$$\|\calA-\widetilde{\calA}\|^2\to\min.$$

\section{CP-ALS with a posteriori antisymmetrization}\label{sec:algorithm1}

First we describe a naive approach. The process is made of two steps. Step $1$: Using the CP-ALS algorithm~\ref{agm:cpals} that ignores the tensor structure we find a rank-$1$ approximation $\bar{\calA}$ of $\calA$,
$$\bar{\calA}=[[x,y,z]], \quad \ran(\bar{\calA})=1.$$
Step $2$: We apply the antisymmetrizer~\eqref{anti3d} on $\bar{\calA}$ to obtain $\widetilde{\calA}$ in the form~\eqref{antiform},
$$\widetilde{\calA}=\anti(\bar{\calA}),$$
that is,
$$\widetilde{\calA}=\calA_6(x,y,z).$$
This procedure is given in Algorithm~\ref{agm:cp+antisym}. We do not need to form the tensor $\bar{\calA}$ explicitly.
\bigskip

\begin{Algorithm}\label{agm:cp+antisym}
\vspace{0.5ex}\hrule\vspace{0.5ex}
\emph{CP with a posteriori antisymmetrization}
\vspace{0.5ex}\hrule
\begin{algorithmic}
\State \textbf{Input:} $\calA\in\R^{n\times n\times n}$ antisymmetric
\State \textbf{Output:} $\widetilde{\calA}=\calA_6(x,y,z)$
\State Apply Algorithm~\ref{agm:cpals} on $\calA$ with $r=1$ to obtain $x,y,z\in\R^n$
\State $\widetilde{\calA}=\calA_6(x,y,z)$
\end{algorithmic}
\hrule
\end{Algorithm}
\bigskip

Obviously, using a rank-$1$ intermediate tensor produces an unnecessarily big approximation error. However, it can be easily shown that if the error of the rank-$1$ approximation is bounded by some $\epsilon > 0$, the resulting error will also be bounded by $\epsilon$.

\section{Antisymmetry-preserving CP algorithm}\label{sec:algorithm2}

For a given antisymmetric tensor $\calA\in\mathbb{R}^{n\times n\times n}$ we are looking for vectors $x,y,z\in\R^{n}$ such that
\begin{equation}\label{problem}
\|\calA-\calA_6(x,y,z)\|^2\to\min.
\end{equation}
Contrary to the Algorithm~\ref{agm:cpals}, here we develop a new structure-preserving low-rank appro\-xi\-ma\-ti\-on algorithm. Our algorithm uses the ALS approach, that is, we are solving an optimization problem in each mode. It results with a tensor of the form~\eqref{antiform} and there is no need to apply the antisymmetrizer.
ALS algorithms are widely used to address different multilinear minimization problems~\cite{CLA09,Usch12,EHK15,GKK15}, including the ones regarding the CP appro\-xi\-ma\-ti\-on~\cite{AcDuKo11,MVLB23,HTG17}. There is also a very recent extension to the antisymmetric case~\cite{TFS23}, but both the problem and the algorithm are different from ours.

Set
$$a=\left[
         \begin{array}{c}
           x \\
           y \\
           z \\
         \end{array}
       \right]\in\R^{3n}.$$
Then, similarly to what was done in~\cite{AcDuKo11}, we define the objective function $f\colon\R^{3n}\to\R$,
\begin{equation}\label{function}
f(a)=6\|\calA-\calA_6(x,y,z)\|^2.
\end{equation}
We consider three partial minimization problems:
\begin{equation}\label{minimization-xyz}
\min_x f(a), \quad \min_y f(a), \quad \min_z f(a).
\end{equation}

Before we formulate the algorithm, we need to prove Theorem~\ref{tm:alg2}. It gives three reformulations of the objective function $f$ that we are going to use in order to find
the solutions of the problems~\eqref{minimization-xyz}.

Observe that, since $\calA_6(x,y,z)$ is linear in $x$, $y$ and $z$, the objective function is quadratic in $x$, $y$ and $z$. The approximation problem becomes a quadratic optimization problem. Here we derive the quadratic forms explicitly. However, it is worth mentioning that the underlying linearity opens the possibilities of the extension to more general settings.

In order to simplify the statement of the theorem, we define the following objects: matrices $Q^{(1)}=Q^{(1)}(y,z),Q^{(2)}=Q^{(2)}(x,z),Q^{(3)}=Q^{(3)}(x,y)\in\R^{n\times n}$,
\begin{align}
Q^{(1)} & = 2\left(\left(\|y\|_2^2\|z\|_2^2-\langle y,z\rangle^2\right)I_n+(yz^T-zy^T)^2\right), \label{tm:Qx} \\
Q^{(2)} & = 2\left(\left(\|z\|_2^2\|x\|_2^2-\langle z,x\rangle^2\right)I_n+(zx^T-xz^T)^2\right), \label{Qy} \\
Q^{(3)} & = 2\left(\left(\|x\|_2^2\|y\|_2^2-\langle x,y\rangle^2\right)I_n+(xy^T-yx^T)^2\right); \label{Qz}
\end{align}
vectors $c^{(1)}=c^{(1)}(y,z),c^{(2)}=c^{(2)}(x,z),c^{(3)}=c^{(3)}(x,y)\in\R^n$, 
\begin{align}
c^{(1)} & = -12\calA\times_2y^T\times_3z^T, \label{tm:cx} \\
c_i^{(2)} & = -12\calA\times_2z^T\times_3x^T, \label{cy} \\
c_i^{(3)} & = -12\calA\times_2x^T\times_3y^T; \label{cz}
\end{align}
and a real number
\begin{equation}\label{tm:d}
d=6\|\calA\|^2.
\end{equation}

\begin{Theorem}\label{tm:alg2}
The function $f$ defined by~\eqref{function} can be written as
\begin{align}
f(a) & = d+(c^{(1)})^Tx +\frac{1}{2}x^TQ^{(1)}x \label{tm:x}\\
& = d+(c^{(2)})^Ty +\frac{1}{2}y^TQ^{(2)}y \label{tm:y}\\
& = d+(c^{(3)})^Tz +\frac{1}{2}z^TQ^{(3)}z, \label{tm:z}
\end{align}
for $Q^{(1)},Q^{(2)},Q^{(3)}\in\R^{n\times n}$, $c^{(1)},c^{(2)},c^{(3)}\in\R^n$, and $d\in\R$ defined by the relations~\eqref{tm:Qx}--\eqref{tm:d}.
\end{Theorem}

\begin{proof}
First, we can write the function $f$ from~\eqref{function} as
\begin{align}
f(a) & =6\|\calA\|^2 -2\langle\calA,6\calA_6(x,y,z)\rangle +\frac{1}{6}\|6\calA_6(x,y,z)\|^2 \nonumber \\
& =6f_1(a) -2f_2(a) +\frac{1}{6}f_3(a), \label{f123}
\end{align}
where
\begin{align}
f_1(a) & =\|\calA\|^2, \nonumber \\
f_2(a) & =\langle\calA,x\circ y\circ z+y\circ z\circ x+z\circ x\circ y-x\circ z\circ y-y\circ x\circ z-z\circ y\circ x\rangle, \label{f2} \\
f_3(a) & =\|x\circ y\circ z+y\circ z\circ x+z\circ x\circ y-x\circ z\circ y-y\circ x\circ z-z\circ y\circ x\|^2. \label{f3}
\end{align}

For the function $f_2$ we have
\begin{align*}
f_2(a) & = \sum_{i,j,k=1}^n a_{ijk}(x_iy_jz_k+y_iz_jx_k+z_ix_jy_k-x_iz_jy_k-y_ix_jz_k-z_iy_jx_k) \\
& = \sum_{i=1}^n x_i \sum_{j,k=1}^n a_{ijk}y_jz_k +\sum_{k=1}^n x_k \sum_{i,j=1}^n a_{ijk}y_iz_j + \sum_{j=1}^n x_j \sum_{i,k=1}^n a_{ijk}z_iy_k \\
& \ +\sum_{i=1}^n x_i \sum_{j,k=1}^n (-a_{ijk})z_jy_k +\sum_{k=1}^n x_k \sum_{i,j=1}^n (-a_{ijk})z_iy_j +\sum_{j=1}^n x_j \sum_{i,k=1}^n (-a_{ijk})y_iz_k.
\end{align*}
We rename the indices in the upper expression and use the fact that $\calA$ is antisymmetric. We get
\begin{align}
f_2(a) & = \sum_{i=1}^n x_i \sum_{j,k=1}^n a_{ijk}y_jz_k +\sum_{i=1}^n x_i \sum_{j,k=1}^n a_{jki}y_jz_k +\sum_{i=1}^n x_i \sum_{j,k=1}^n a_{kij}y_jz_k \nonumber \\
& \ +\sum_{i=1}^n x_i \sum_{j,k=1}^n (-a_{ikj})y_jz_k +\sum_{i=1}^n x_i \sum_{j,k=1}^n (-a_{kji})y_jz_k +\sum_{i=1}^n x_i \sum_{j,k=1}^n (-a_{jik})y_jz_k \nonumber \\
& = 6\sum_{i=1}^n x_i \sum_{j,k=1}^n a_{ijk}y_jz_k. \label{f2x}
\end{align}

Further on, we write the function $f_3$ as
\begin{equation}\label{f3a}
f_3(a) = \sum_{i,j,k=1}^n(x_iy_jz_k+y_iz_jx_k+z_ix_jy_k-x_iz_jy_k-y_ix_jz_k-z_iy_jx_k)^2.
\end{equation}
After regrouping the summands and renaming the indices, as we did for $f_2$, it follows from~\eqref{f3a} that
\begin{align*}
& f_3(a) = \\
& = 6\sum_{i=1}^n x_i^2\left(\sum_{j,k=1}^n y_j^2z_k^2\right) - 6\sum_{i=1}^n x_i^2\left(\sum_{j,k=1}^n y_jy_kz_jz_k\right) \\
& \ +12\sum_{i,j=1}^n x_ix_j\left(\sum_{k=1}^n y_iy_kz_jz_k\right) - 6\sum_{i,j=1}^n x_ix_j\left(\sum_{k=1}^n y_iy_jz_k^2\right)
- 6\sum_{i,j=1}^n x_ix_j\left(\sum_{k=1}^n y_k^2z_iz_j\right) \\
& = \sum_{i=1}^n x_i^2 \left(6\sum_{j=1}^n y_j^2\sum_{k=1}^n z_k^2 - 6\bigg(\sum_{j=1}^n y_jz_j\bigg)^2\right) \\
& \ + \sum_{i,j=1}^n x_ix_j \left(12\sum_{k=1}^n y_iy_kz_jz_k - 6\sum_{k=1}^n y_iy_jz_k^2 - 6\sum_{k=1}^n y_k^2z_iz_j\right) \\
& = \sum_{i=1}^n x_i^2 \left(6\sum_{j=1}^n y_j^2\sum_{k=1}^n z_k^2 - 6\bigg(\sum_{j=1}^n y_jz_j\bigg)^2 + 12y_iz_i\sum_{k=1}^n y_kz_k - 6y_i^2\sum_{k=1}^n z_k^2 - 6z_i^2\sum_{k=1}^n y_k^2\right) \\
& \ + \sum_{\substack{i,j=1 \\ i<j}}^n x_ix_j \left(12y_iz_j\sum_{k=1}^n y_kz_k + 12z_iy_j\sum_{k=1}^n y_kz_k - 12y_iy_j\sum_{k=1}^n z_k^2 - 12z_iz_j\sum_{k=1}^n y_k^2\right).
\end{align*}
That is,
\begin{align}
f_3(a) & = \sum_{i=1}^n x_i^2 \left(6\|y\|_2^2\|z\|_2^2 - 6\langle y,z\rangle^2 + 12y_iz_i\langle y,z\rangle - 6y_i^2\|z\|_2^2 - 6z_i^2\|y\|_2^2\right) \nonumber \\
& \ + \sum_{\substack{i,j=1 \\ i<j}}^n x_ix_j \left(12(y_iz_j+z_iy_j)\langle y,z\rangle - 12y_iy_j\|z\|_2^2 - 12z_iz_j\|y\|_2^2\right). \label{f3x}
\end{align}

Then, we can set
\begin{align*}
d & = 6f_1(a), \\
(c^{(1)})^Tx & = -2f_2(a), \\
\frac{1}{2}x^TQ^{(1)}x & = \frac{1}{6}f_3(a).
\end{align*}
From the relations~\eqref{f123}, \eqref{f2x}, and~\eqref{f3x} we get the assertion~\eqref{tm:x} where
\begin{equation}\label{cx}
c_i^{(1)} = -12\sum_{j,k=1}^n a_{ijk}y_jz_k, \quad 1\leq i\leq n,
\end{equation}
\begin{align}
q_{ii}^{(1)} & = 2\|y\|_2^2\|z\|_2^2 - 2\langle y,z\rangle^2 + 4y_iz_i\langle y,z\rangle - 2y_i^2\|z\|_2^2 - 2z_i^2\|y\|_2^2, \nonumber \\
q_{ij}^{(1)} & = 2(y_iz_j+z_iy_j)\langle y,z\rangle - 2y_iy_j\|z\|_2^2 - 2z_iz_j\|y\|_2^2, \quad 1\leq i,j\leq n, \ i\neq j, \label{Qx}
\end{align}
and $d$ is as given in~\eqref{tm:d}.
It follows from the expressions in~\eqref{Qx} that
\begin{align*}
Q^{(1)} & = 2\left(\|y\|_2^2\|z\|_2^2-\langle y,z\rangle^2\right)I_n
+ 2\left((yz^T+zy^T)\langle y,z\rangle - yy^T\|z\|_2^2-zz^T\|y\|_2^2\right) \\
& = 2\left(\left(\|y\|_2^2\|z\|_2^2-\langle y,z\rangle^2\right)I_n+(yz^T-zy^T)^2\right),
\end{align*}
while the vector given element-wise by~\eqref{cx} is equal to the one from the relation~\eqref{tm:cx}.

Similarly, using a different regrouping of the summands in the equations~\eqref{f2} and~\eqref{f3}, we obtain the assertions~\eqref{tm:y} and~\eqref{tm:z}. We get $Q^{(2)}$ and $c^{(2)}$, as in the relations~\eqref{Qy} and~\eqref{cy}, respectively, as well as $Q^{(3)}$ and $c^{(3)}$, as in~\eqref{Qz} and~\eqref{cz}, respectively.
\end{proof}

Minimization problem of the form
$$\min_v\{d+c^Tv +\frac{1}{2}v^TQv\},$$
is a problem of quadratic programming with no constraints.
Its solution $v$ is given by the linear system
$$Qv=-c.$$
Therefore, in order to find the solutions of the minimization problems
\begin{equation}\label{min-xyz}
\left.
  \begin{array}{c}
    \displaystyle \min_x d+(c^{(1)})^Tx +\frac{1}{2}x^TQ^{(1)}x \\
    \displaystyle \min_y d+(c^{(2)})^Ty +\frac{1}{2}y^TQ^{(2)}y \\
    \displaystyle \min_z d+(c^{(3)})^Tz +\frac{1}{2}z^TQ^{(3)}z \\
  \end{array}
\right\}
\end{equation}
we need to solve linear systems
\begin{equation*}
\left.
  \begin{array}{c}
    Q^{(1)}x=-c^{(1)} \\
    Q^{(2)}y=-c^{(2)} \\
    Q^{(3)}z=-c^{(3)} \\
  \end{array}
\right\}
\end{equation*}
respectively.

Here we come to an obstacle because matrices $Q^{(1)},Q^{(2)},Q^{(3)}$ are singular. Take $Q^{(1)}$. From the relation~\eqref{tm:Qx} we see that $Q^{(1)}$ is defined by two vectors $y$ and $z$ and we have
$$Q^{(1)}y=0, \quad Q^{(1)}z=0.$$
Assuming that $y$ and $z$ are linearly independent vectors, this means that $\ran(Q^{(1)})\leq n-2$. On the other hand, $Q^{(1)}$ is defined as an identity matrix minus a rank-$2$ matrix. This implies that $\ran(Q^{(1)})=n-2$. However, linear system $Q^{(1)}x=-c^{(1)}$ is consistent because
$\ran([Q^{(1)}c^{(1)}])=\ran(Q^{(1)})$, which can be seen from the relations~\eqref{tm:Qx} and~\eqref{tm:cx}. Hence, linear system $Q^{(1)}x=-c^{(1)}$ can be solved using the Moore-Penrose inverse,
$$x=-\left(Q^{(1)}\right)^{+}c^{(1)},$$
with an additional constraint that $x$ is orthogonal to $y$ and $z$.
Analogue reasoning holds for the linear systems for $y$ and $z$.

Now, we can write the algorithm for solving the minimization problem~\eqref{problem}. The algorithm is based on solving three minimization problems~\eqref{min-xyz}.
\bigskip

\begin{Algorithm}\label{agm:antisymcp}
\vspace{0.5ex}\hrule\vspace{0.5ex}
\emph{Antisymmetry-preserving CP}
\vspace{0.5ex}\hrule
\begin{algorithmic}
\State \textbf{Input:} $\calA\in\R^{n\times n\times n}$ antisymmetric
\State \textbf{Output:} $\widetilde{\calA}=\calA_6(x,y,z)$
\State Initialize $x,y,z\in\R^n$ as random vectors.
\Repeat
\State For $c^{(1)}$ as in~\eqref{tm:cx} and $Q^{(1)}$ as in~\eqref{tm:Qx},
$x=-\left(Q^{(1)}\right)^{+}c^{(1)}$.
\State For $c^{(2)}$ as in~\eqref{cy} and $Q^{(2)}$ as in~\eqref{Qy},
$y=-\left(Q^{(2)}\right)^{+}c^{(2)}$.
\State For $c^{(3)}$ as in~\eqref{cz} and $Q^{(3)}$ as in~\eqref{Qz},
$z=-\left(Q^{(3)}\right)^{+}c^{(3)}$.
\Until{convergence or maximum number of iterations}
\State $\widetilde{\calA}=\calA_6(x,y,z)$
\end{algorithmic}
\hrule
\end{Algorithm}
\bigskip

Note that the Algorithm~\ref{agm:antisymcp} results with mutually orthogonal vectors $x$, $y$ and $z$. While this may seem restrictive, it is reasonable to do a restriction on orthogonal vectors. Let us explain that.

First, note that if $x$, $y$ and $z$ are linearly dependent, then $\calA_6(x,y,z)=0$. Thus, any linearly independent triplet of vectors will give a smaller value of the objective function $f$ defined by the relation~\eqref{function}.

The objective function $f$ is invariant under very general transformations.
Due to multilinearity, for $\alpha,\beta\in\R$, we have
$$\calA_6(x+\alpha y+\beta z,y,z)=\calA_6(\alpha y+\beta z,y,z)+\calA_6(x,y,z)=\calA_6(x,y,z),$$
where $\calA_6(\alpha y+\beta z,y,z)=0$ by antisymmetry.
Let us consider an arbitrary matrix $B=(b_{ij})\in\R^{3\times3}$. Using the arguments of multilinearity and antisymmetry as above, we have
\begin{align*}
& \calA_6(b_{11}x+b_{12}y+b_{13}z,b_{21}x+b_{22}y+b_{23}z,b_{31}x+b_{32}y+b_{33}z) = \\
& = b_{11}b_{22}b_{33}\calA_6(x,y,z)+b_{11}b_{23}b_{32}\calA_6(x,z,y)+b_{12}b_{21}b_{33}\calA_6(y,x,z) + \\
& \quad +b_{12}b_{23}b_{31}\calA_6(y,z,x)+b_{13}b_{21}b_{32}\calA_6(z,x,y)+b_{13}b_{22}b_{31}\calA_6(z,y,x) \\
& = (b_{11}b_{22}b_{33}-b_{11}b_{23}b_{32}-b_{12}b_{21}b_{33}+b_{12}b_{23}b_{31}+b_{13}b_{21}b_{32}-b_{13}b_{22}b_{31})\calA_6(x,y,z) \\
& = \de(B)\calA_6(x,y,z).
\end{align*}
Therefore, if $\de(B)=1$,
$$\|\calA-\calA_6(b_{11}x+b_{12}y+b_{13}z,b_{21}x+b_{22}y+b_{23}z,b_{31}x+b_{32}y+b_{33}z)\|=
\|\calA-\calA_6(x,y,z)\|,$$
and the value of the objective function $f$ stays the same.

Set $V=[x,y,z]\in\R^{n\times3}$. Take the thin QR decomposition $V=\widetilde{V}R$, such that $det(R)=1$ and $\widetilde{V}=[\tilde{x},\tilde{y},\tilde{z}]\in\R^{n\times3}$ has orthogonal columns. Then, following the same reasoning, we have
$$\calA_6(\tilde{x},\tilde{y},\tilde{z})=\calA_6(x,y,z)$$
and
\begin{equation}\label{min-orth}
\min_{x,y,z}\|\calA-\calA_6(x,y,z)\|=\min_{\tilde{x},\tilde{y},\tilde{z} \text{ orthogonal}}\|\calA-\calA_6(\tilde{x},\tilde{y},\tilde{z})\|.
\end{equation}
This justifies the choice of orthogonal vectors.

\subsection{Equivalence of the Algorithm~\ref{agm:antisymcp} to HOPM}\label{sec:hopm}

Here we are going to show that the Algorithm~\ref{agm:antisymcp} is equivalent to the higher-order power method (HOPM) for unstructured rank-$1$ approximation and see what implications it has.

Due to multilinearity, minimization problem~\eqref{min-orth} can be modified into a minimization problem on unitary vectors,
\begin{equation}\label{min-orth2}
\min_{\substack{\tilde{x},\tilde{y},\tilde{z} \text{ orthonormal} \\ \lambda\in\R}}\|\calA-\lambda\calA_6(\tilde{x},\tilde{y},\tilde{z})\|^2,
\end{equation}
so it becomes a minimization problem on the Stiefel manifold. Since the expression in~\eqref{min-orth2} does not depend on the basis, it is a minimization problem on the Grassmann manifold, which makes sense as the antisymmetric tensors are connected to the Grassmannians, see, e.g.,~\cite{Lan12}.

We can rewrite~\eqref{min-orth2} as
\begin{equation}\label{min-orth3}
\min_{\substack{\tilde{x},\tilde{y},\tilde{z} \text{ orthonormal} \\ \lambda\in\R}}
\left\{\|\calA\|^2-2\lambda\langle\calA,\calA_6(\tilde{x},\tilde{y},\tilde{z})\rangle+\lambda^2\|\calA_6(\tilde{x},\tilde{y},\tilde{z})\|^2\right\}.
\end{equation}
Set
\begin{equation}\label{matrixV}
V=\left[
  \begin{array}{ccc}
    \tilde{x} & \tilde{y} & \tilde{z} \\
  \end{array}
\right].
\end{equation}
Observe that $$\calA_6(\tilde{x},\tilde{y},\tilde{z})=\calE\times_1V\times_2V\times_3V,$$
where $\calE$ is given by the relation~\eqref{tensorE}. Then,
$$\|\calA_6(\tilde{x},\tilde{y},\tilde{z})\|^2=\|\calE\|^2=6,$$
because $\tilde{x}$, $\tilde{y}$, $\tilde{z}$ are orthonormal and the Frobenius norm is unitary invariant.
This way, minimization problem~\eqref{min-orth3} is simplified to
$$\min_{\substack{\tilde{x},\tilde{y},\tilde{z} \text{ orthonormal} \\ \lambda\in\R}}
\left\{\|\calA\|^2-2\lambda\langle\calA,\calA_6(\tilde{x},\tilde{y},\tilde{z})\rangle+6\lambda^2\right\}.$$

Take the objective function \begin{equation}\label{min-orthf}\tilde{f}(\lambda,\tilde{x},\tilde{y},\tilde{z})=\|\calA\|^2-2\lambda\langle\calA,\calA_6(\tilde{x},\tilde{y},\tilde{z})\rangle+6\lambda^2.
\end{equation}
In order to find the optimal $\lambda_*$ for $\tilde{f}$ we set the partial derivative of $\tilde{f}$ to zero. We have
$$\frac{\partial}{\partial\lambda}\tilde{f}(\lambda,\tilde{x},\tilde{y},\tilde{z})=12\lambda-2\langle\calA,\calA_6(\tilde{x},\tilde{y},\tilde{z})\rangle=0,$$
that is,
$$\lambda_*=\frac{\langle\calA,\calA_6(\tilde{x},\tilde{y},\tilde{z})\rangle}{6}.$$
It follows from~\eqref{min-orthf} that
$$\tilde{f}(\lambda_*,\tilde{x},\tilde{y},\tilde{z})=\|\calA\|^2-\frac{1}{6}\langle\calA,\calA_6(\tilde{x},\tilde{y},\tilde{z})\rangle^2.$$
Thus, minimizing $\tilde{f}(\lambda_*,\tilde{x},\tilde{y},\tilde{z})$ is equivalent to maximizing $|\langle\calA,\calA_6(\tilde{x},\tilde{y},\tilde{z})\rangle|$ over the Stiefel manifold.

Define the compressed tensor
\begin{equation}\label{Ac}
\calA_c(V)\coloneqq\calA\times_1V^T\times_2V^T\times_3V^T,
\end{equation}
where $V$ is as in the relation~\eqref{matrixV}. This is a $3\times3\times3$ tensor. It is very similar to tensor $\calE$, except that in place of $1$ and $-1$ it has $(\calA_c(V))_{123}$ and $-(\calA_c(V))_{123}$, respectively.
Using this tensor we obtain
\begin{align*}
|\langle\calA,\calA_6(\tilde{x},\tilde{y},\tilde{z})\rangle| & =|\langle\calA,\calE\times_1V\times_2V\times_3V\rangle|=|\langle\calA_c(V),\calE\rangle| \\
& =6|(\calA_c(V))_{123}| = \sqrt{6}\|\calA_c(V)\|.
\end{align*}
In the last equation we used the norm of the compressed tensor, $\|\calA_c(V)\|^2=6((\calA_c(V))_{123})^2$.
Therefore, maximization of $|\langle\calA,\calA_6(\tilde{x},\tilde{y},\tilde{z})\rangle|$ is equivalent to maximization of $\|\calA_c(V)\|$.
This corresponds to the best structure-preserving multilinear rank-$r$ approximation from~\cite{BeKre17} for $r=3$.

The problem of finding the best antisymmetric multilinear rank-$r$ approximation is equivalent to the problem of finding the best unstructured rank-$1$ approximation of an antisymmetric tensor, see~\cite[Theorem 4.2]{BeKre17}. This implies the equivalence between our Algorithm~\ref{agm:antisymcp} and HOPM used for finding the best unstructured rank-$1$ approximation. Finally, the global convergence result for HOPM given in~\cite{Usch15}, namely, the iterates of the ALS algorithm for HOPM converge to the stationary point of the corresponding objective function, apply to our algorithm as well.

\section{Partial antisymmetry}\label{sec:partial}

Regarding the antisymmetric tensors, we can ask what happens if a tensor has only partial antisymmetry.
We observe order-$3$ tensors. Note that partially antisymmetric tensors do not need to be cubical.

Tensor $\calC\in\R^{n\times n\times m}$ is antisymmetric in modes one and two if all its frontal slices are antisymmetric. Without the loss of generality, we assume that tensor $\calC$ is antisymmetric in the first two modes. That is,
\begin{equation}\label{partial}
c_{ijk}=-c_{jik}, \quad 1\leq i,j\leq n, \ 1\leq k\leq m.
\end{equation}
Tensors that are antisymmetric in modes two and three, or in modes one and three, are defined correspondingly.
Partial antisymmetrizer that results with the antisymmetry in modes one and two can be defined as the operator $\anti_{1,2}$ such that for $\calB\in\R^{n\times n\times m}$ and $\calC=\anti_{1,2}(\calB)$ we have
$$c_{ijk}=\frac{1}{2}(b_{ijk}-b_{jik}).$$

For a pair of indices $(i,j)$, $1\leq i<j\leq n$, subtensor $\calG$ of $\calC$ obtained at the intersection of the $i$th and $j$th column, row, and tube is a $2\times2\times2$ tensor of the form
$$\calG(i_1,i_2,i_3)=\left\{
                       \begin{array}{rl}
                         \alpha, & \text{if $(i_1,i_2,i_3)=(1,2,1)$,} \\
                         -\alpha, & \text{if $(i_1,i_2,i_3)=(2,1,1)$,} \\
                         \beta, & \text{if $(i_1,i_2,i_3)=(1,2,2)$,} \\
                         -\beta, & \text{if $(i_1,i_2,i_3)=(2,1,2)$,} \\
                         0, & \text{if $i_1=i_2$,} \\
                       \end{array}
                     \right.$$
for $\alpha,\beta\in\R$.
Its mode-$1$ matricization is given by
$$G_{(1)}=\left[
                  \begin{array}{cc|cc}
                    0 & \alpha & 0 & \beta \\
                    -\alpha & 0 & -\beta & 0 \\
                  \end{array}
\right].
$$
Here, tensor $\calG$ plays the role analogue to the Levi-Civita tensor~\eqref{tensorE} in Section~\ref{sec:problem}.

Analogue to the tensor format~\eqref{antiform}, for three vectors $x,y\in\R^{n}$ and $z\in\R^{m}$, we can define an $n\times n\times m$ tensor
\begin{equation}\label{antiform_partial}
\calC_2(x,y,z):=\frac{1}{2}(x\circ y\circ z-y\circ x\circ z).
\end{equation}
If we take $x=[1,0]^T$, $y=[0,1]^T$, $z=[\alpha,\beta]^T$, then $\calC_2(x,y,z)=\calG$.
Besides, if $\calT=[[x,y,z]]$ is a rank-$1$ tensor, then $\calC_2(x,y,z)=\anti_{1,2}(\calT)$.
Obviously, $\ran(\calC_2(x,y,z))\leq2$. For the fixed third index, each slice of $\calC_2(x,y,z)$ is a skew-symmetric matrix and, therefore, has an even rank. Hence, $$\ran(\calC_2(x,y,z))=2.$$

Considering all this, for a given non-zero tensor $\calC\in\mathbb{R}^{n\times n\times m}$ that is antisymmetric in the first two modes, we are looking for its rank-$2$ approximation $\widetilde{\calC}$ of the same structure.
Again, we examine two approaches. The first one is analogue to the Section~\ref{sec:algorithm1}.
In the second approach we find a tensor $\widetilde{\calC}=\calC_2(x,y,z)$ defined by the vectors $x,y\in\R^{n}$ and $z\in\R^{m}$,
such that
\begin{equation}\label{minimization_partial}
\|\calC-\widetilde{\calC}\|^2\to\min.
\end{equation}

\subsection{Ignoring the structure}

Let $\calC\in\R^{n\times n\times m}$ ba a tensor with partial antisymmetry. We first approximate $\calC$ with a rank-$1$ tensor $\bar{\calC}$ by using the CP-ALS algorithm~\ref{agm:cpals} with $r=1$. Then, we apply the operator $\anti_{1,2}$ on $\bar{\calC}$ to get a rank-$2$ tensor $\widetilde{\calC}$ that is antisymmetric in modes one and two. We have
\begin{align*}
\bar{\calC} & =[[x,y,z]], \quad x,y\in\R^n, \ z\in\R^m, \\
\widetilde{\calC} & = \anti_{1,2}(\bar{\calC}),
\end{align*}
or equivalently, $\widetilde{\calC}=\calC_2(x,y,z)$.
The algorithm with partial a posteriori antisymmetrization is a simple modification of the Algorithm~\ref{agm:cp+antisym}.
\bigskip

\begin{Algorithm}\label{agm:cp+pantisym}
\vspace{0.5ex}\hrule\vspace{0.5ex}
\emph{CP with partial a posteriori antisymmetrization}
\vspace{0.5ex}\hrule
\begin{algorithmic}
\State \textbf{Input:} $\calC\in\R^{n\times n\times m}$ antisymmetric in modes $1$ and $2$
\State \textbf{Output:} $\widetilde{\calC}=\calC_2(x,y,z)$
\State Apply Algorithm~\ref{agm:cpals} on $\calA$ with $r=1$ to obtain $x,y\in\R^n$, $z\in\R^m$
\State $\widetilde{\calC}=\calC_2(x,y,z)$
\end{algorithmic}
\hrule
\end{Algorithm}
\bigskip

\subsection{Preserving the structure}

Now we are going to construct an iterative structure-preserving minimization algorithm. Again, let $\calC\in\R^{n\times n\times m}$ be a tensor with partial antisymmetry. We are looking for tensor $\widetilde{\calC}\in\R^{n\times n\times m}$ which is a solution of the minimization problem~\eqref{minimization_partial}. In particular, we are looking for vectors $x,y\in\R^n$ and $z\in\R^m$ such that
\begin{equation}\label{problem_partial}
\|\calC-\calC_2(x,y,z)\|^2\to\min.
\end{equation}

We set
$$v=\left[
      \begin{array}{c}
        x \\
        y \\
        z \\
      \end{array}
    \right]\in\R^{2n+m},
$$
and define the objective function $g\colon\R^{2n+m}\to\R$,
\begin{equation}\label{function_partial}
g(v)=2\|\calC-\calC_2(x,y,z)\|^2.
\end{equation}
We formulate the ALS algorithm based on three minimization problems:
$$\min_x g(v), \quad \min_y g(v), \quad \min_z g(v).$$

To this end, we need Theorem~\ref{tm:alg4}. Before the statement of the theorem, we define the appropriate objects: matrices $Q^{(1)}=Q^{(1)}(y,z),Q^{(2)}=Q^{(2)}(x,z)\in\R^{n\times n}$,
\begin{align}
Q^{(1)} & = 2\|y\|_2^2\|z\|_2^2I_n-2yy^T\|z\|_2^2,\label{tm:Qxp} \\
Q^{(2)} & = 2\|x\|_2^2\|z\|_2^2I_n-2xx^T\|z\|_2^2; \label{Qyp}
\end{align}
vectors $b^{(1)}=b^{(1)}(y,z),b^{(2)}=b^{(2)}(x,z)\in\R^n$, $b^{(3)}=b^{(3)}(x,y)\in\R^m$,
\begin{align}
b^{(1)} & = -4\calC\times_2y^T\times_3z^T, \label{tm:cxp} \\
b^{(2)} & = -4\calC\times_2x^T\times_3z^T, \label{cyp} \\
b^{(3)} & = -2(\calC\times_1x^T\times_2y^T-\calC\times_1y^T\times_2x^T), \label{tm:czp}
\end{align}
and numbers $q^{(3)}=q^{(3)}(x,y),d\in\R$,
\begin{align}
q^{(3)} & = \|xy^T-yx^T\|_2^2, \label{Qzp} \\
d & = 2\|\calC\|^2. \label{tm:dp}
\end{align}

\begin{Theorem}\label{tm:alg4}
The function $g$ defined by~\eqref{function_partial} can be written as
\begin{align}
g(v) & = d+(b^{(1)})^Tx +\frac{1}{2}x^TQ^{(1)}x \label{tm:xp}\\
& = d+(b^{(2)})^Ty +\frac{1}{2}y^TQ^{(2)}y \label{tm:yp}\\
& = d+(b^{(3)})^Tz +\frac{1}{2}q^{(3)}z^Tz, \label{tm:zp}
\end{align}
for $Q^{(1)},Q^{(2)}\in\R^{n\times n}$, $b^{(1)},b^{(2)}\in\R^n$, $b^{(3)}\in\R^m$, $q^{(3)}\in\R$ defined by the relations~\eqref{tm:Qxp}-\eqref{tm:dp}.
\end{Theorem}

\begin{proof}
We start by writing the function $g$ as
\begin{align*}
g(v) & =2\|\calC\|^2-2\left\langle \calC,x\circ y\circ z-y\circ x\circ z\right\rangle +
\frac{1}{2}\|x\circ y\circ z-y\circ x\circ z\|^2 \\
& =2g_1(v)-2g_2(v)+\frac{1}{2}g_3(v),
\end{align*}
for
\begin{align}
g_1(v) & = \|\calC\|^2, \nonumber \\
g_2(v) & = \left\langle \calC,x\circ y\circ z-y\circ x\circ z\right\rangle, \label{g2}\\
g_3(v) & = \|x\circ y\circ z-y\circ x\circ z\|^2. \label{g3}
\end{align}

Function $g_2$ can be written as
\begin{align*}
g_2(v) & = \sum_{i,j=1}^n\sum_{k=1}^m c_{ijk}(x_iy_jz_k-y_ix_jz_k) \\
& = \sum_{i=1}^n x_i\left(\sum_{j=1}^n\sum_{k=1}^m c_{ijk}y_jz_k\right) + \sum_{j=1}^n x_j\left(\sum_{i=1}^n\sum_{k=1}^m (-c_{ijk})y_iz_k\right).
\end{align*}
Using the partial antisymmetry property~\eqref{partial}, after renaming the indices we get
$$g_2(v) = 2\sum_{i=1}^n x_i\left(\sum_{j=1}^n\sum_{k=1}^m c_{ijk}y_jz_k\right).$$
For the function $g_3$ we have
\begin{align*}
g_3(v) & = \sum_{i,j=1}^n\sum_{k=1}^m (x_iy_jz_k-x_jy_iz_k)^2 \\
& = \sum_{i=1}^n x_i^2\left(\sum_{j=1}^n\sum_{k=1}^m y_j^2z_k^2\right) - 2\sum_{i,j=1}^n x_ix_jy_iy_j\left(\sum_{k=1}^m z_k^2\right) +
\sum_{j=1}^n x_j^2\left(\sum_{i=1}^n\sum_{k=1}^m y_i^2z_k^2\right) \\
& = 2\sum_{i=1}^n x_i^2\|y\|_2^2\|z\|_2^2 - 2\sum_{i,j=1}^n x_ix_jy_iy_j\|z\|_2^2 \\
& = \sum_{i=1}^n x_i^2(2\|y\|_2^2\|z\|_2^2-2y_i^2\|z\|_2^2) + \sum_{\substack{i,j=1 \\ i\neq j}}^n x_ix_j(-2y_iy_j\|z\|_2^2).
\end{align*}

The same way as in the proof of Theorem~\ref{tm:alg2}, we set
\begin{align*}
d & = 2g_1(v), \\
(b^{(1)})^Tx & = -2g_2(v), \\
\frac{1}{2}x^TQ^{(1)}x & = \frac{1}{2}g_3(v),
\end{align*}
where
\begin{equation}\label{cxp}
b_i^{(1)} = -4\sum_{j=1}^n\sum_{k=1}^m c_{ijk}y_jz_k, \quad 1\leq i\leq n,
\end{equation}
and
\begin{align}
q_{ii}^{(1)} & = 2\|y\|^2\|z\|^2-2y_i^2\|z\|^2, \nonumber \\
q_{ij}^{(1)} & = -2y_iy_j\|z\|^2, \quad 1\leq i,j\leq n, \ i\neq j. \label{Qxp}
\end{align}
Vector $b^{(1)}$ from~\eqref{cxp} can be written in the more compact form~\eqref{tm:cxp} and matrix $Q^{(1)}$ from~\eqref{Qxp} is equivalent to~\eqref{tm:Qxp}, while $d$ is like in the relation~\eqref{tm:dp}. This is how we get the assertion~\eqref{tm:xp}.

With different regrouping of the summands in the relations~\eqref{g2} and~\eqref{g3} we get the equation~\eqref{tm:yp} with
$b^{(2)}$ and $Q^{(2)}$ as in~\eqref{cyp} and~\eqref{Qyp}, respectively.

To get the equation~\eqref{tm:zp} we write
$$g_2(v) = \sum_{k=1}^m z_k\left(\sum_{i,j=1}^n c_{ijk}(x_iy_j-y_ix_j)\right)$$
and
$$g_3(v) = \sum_{k=1}^m z_k^2\left(\sum_{i,j=1}^n(x_iy_j-x_jy_i)^2\right).$$
Then, we set
\begin{align*}
b_k^{(3)} & = -2\sum_{i,j=1}^n c_{ijk}(x_iy_j-y_ix_j), \quad 1\leq k\leq m, \\
q^{(3)} & = \sum_{i,j=1}^n(x_iy_j-x_jy_i)^2=\|xy^T-yx^T\|_2^2.
\end{align*}
Compact form of the vector $b^{(3)}$ corresponds to~\eqref{tm:czp}.
\end{proof}

Therefore, like in the Section~\ref{sec:algorithm2}, our algorithm is based on finding the solutions of the minimization problems
$$\left.  \begin{array}{c}
    \displaystyle \min_x d+(b^{(1)})^Tx +\frac{1}{2}x^TQ^{(1)}x \\
    \displaystyle \min_y d+(b^{(2)})^Ty +\frac{1}{2}y^TQ^{(2)}y \\
    \displaystyle \min_z d+(b^{(3)})^Tz +\frac{1}{2}q^{(3)}z^Tz. \\
  \end{array}
\right\}$$
Those solutions are obtained, respectively,from the following equations
\begin{equation*}
\left.
  \begin{array}{c}
    Q^{(1)}x=-b^{(1)} \\
    Q^{(2)}y=-b^{(2)} \\
    z=-\frac{1}{q^{(3)}}b^{(3)}. \\
  \end{array}
\right\}
\end{equation*}

The situation regarding these linear systems is similar as for the fully antisymmetric case. Matrices $Q^{(1)}$ and $Q^{(2)}$ are not of the full rank. From their definitions~\eqref{tm:Qxp} and~\eqref{Qyp} we see that both are given as identity minus a rank-$1$ matrix and
$$Q^{(1)}x=0, \quad Q^{(2)}y=0.$$
Thus, $\ran(Q^{(1)})=\ran(Q^{(2)})=n-1.$ Still, $\ran([Q^{(1)}b^{(1)}])=\ran(Q^{(1)})$ and $\ran([Q^{(2)},b^{(2)}])=\ran(Q^{(2)})$, so the linear systems are consistent and can be solved using the Moore-Penrose inverse. Additionally, we get that the vectors $x$ and $y$ must be orthogonal. Note that, for $x\neq y$, $q^{(3)}\neq0$ and $z$ is well-defined.

The algorithm for solving the minimization problem~\eqref{problem_partial} is very similar to the Algorithm~\ref{agm:antisymcp}.
\bigskip

\begin{Algorithm}\label{agm:pantisymcp}
\vspace{0.5ex}\hrule\vspace{0.5ex}
\emph{CP preserving partial antisymmetry}
\vspace{0.5ex}\hrule
\begin{algorithmic}
\State \textbf{Input:} $\calC\in\R^{n\times n\times m}$ antisymmetric in modes $1$ and $2$
\State \textbf{Output:} $\widetilde{\calC}=\calC_2(x,y,z)$
\State Initialize $x,y\in\R^n$,$z\in\R^m$ as random vectors.
\Repeat
\State For $b^{(1)}$ as in~\eqref{tm:cxp} and $Q^{(1)}$ as in~\eqref{tm:Qxp}, $x=-(Q^{(1)})^+b^{(1)}$.
\State For $b^{(2)}$ as in~\eqref{cyp} and $Q^{(2)}$ as in~\eqref{Qyp}, $y=-(Q^{(2)})^+b^{(2)}$.
\State For $b^{(3)}$ as in~\eqref{tm:czp} and $q^{(3)}$ as in~\eqref{Qzp}, $z=-\frac{b^{(3)}}{q^{(3)}}$.
\Until{convergence or maximum number of iterations}
\State $\widetilde{\calC}=\calC_2(x,y,z)$
\end{algorithmic}
\hrule
\end{Algorithm}
\bigskip

Like in the fully antisymmetric case, we can additionally observe that $\calC_2(x,y,z)=0$ if $x$ and $y$ are lineary dependent and
$$\calC_2(b_{11}x+b_{12}y,b_{21}x+b_{22}y,z)=\det B\calC_2(x,y,z), \quad B=\left[
  \begin{array}{cc}
    b_{11} & b_{12} \\
    b_{21} & b_{22} \\
  \end{array}
\right].$$
Then we can rescale our optimization problem such that we are looking for
$$\min_{\substack{\|\tilde{x}\|=\|\tilde{y}\|=\|\tilde{z}\|=1 \\ x\perp y, \ \lambda\in\R}}\left\{\|\calC-\lambda\calC_2(\tilde{x},\tilde{y},\tilde{z})\|^2\right\}.$$
Set
$$\tilde{g}(\lambda,\tilde{x},\tilde{y},\tilde{z})=\|\calC\|^2-2\lambda\langle\calC,\calC_2(\tilde{x},\tilde{y},\tilde{z})\rangle+\lambda^2\|\calC_2(\tilde{x},\tilde{y},\tilde{z})\|^2.$$
From the shape of $\calC_2$ and the fact that $\|\tilde{x}\|=\|\tilde{y}\|=\|\tilde{z}\|=1$ and $x\perp y$, after a short a calculation we get $\|\calC_2(\tilde{x},\tilde{y},\tilde{z})\|^2=\frac{1}{2}$. Thus,
$$\tilde{g}(\lambda,\tilde{x},\tilde{y},\tilde{z})=\|\calC\|^2-2\lambda\langle\calC,\calC_2(\tilde{x},\tilde{y},\tilde{z})\rangle+\frac{1}{2}\lambda^2.$$
The optimal $\lambda$ for $\tilde{g}$ is
$$\lambda_*=2\langle\calC,\calC_2(\tilde{x},\tilde{y},\tilde{z})\rangle$$
and
$$\tilde{g}(\lambda_*,\tilde{x},\tilde{y},\tilde{z}))=\|\calC\|^2-2\langle\calC,\calC_2(\tilde{x},\tilde{y},\tilde{z})\rangle^2.$$
Therefore, minimizing $\tilde{g}(\lambda_*,\tilde{x},\tilde{y},\tilde{z}))$ is equivalent to maximizing $|\langle\calC,\calC_2(\tilde{x},\tilde{y},\tilde{z})\rangle|$.

Now we can set
$$W=\left[
  \begin{array}{cc}
    \tilde{x} & \tilde{y} \\
  \end{array}
\right]$$
and define the compressed matrix
\begin{equation}\label{Cc}
C_c(W,\tilde{z})\coloneqq\calC\times_1W^T\times_2W^T\times_3\tilde{z}^T,
\end{equation}
analogue to the compressed tensor~\eqref{Ac}. Matrix $C_c(W,\tilde{z})$ is a $2\times2$ skew-symmetric matrix
$$\left[
  \begin{array}{cc}
    0 & (C_c(W,\tilde{z}))_{12} \\
    -(C_c(W,\tilde{z}))_{12} & 0 \\
  \end{array}
\right],$$
where
\begin{equation}\label{Ccvalue}
|(C_c(W,\tilde{z}))_{12}|=|\calC\times_1\tilde{x}^T\times_2\tilde{y}^T\times_3\tilde{z}^T|.
\end{equation}
Moreover, we can write
$$\calC_2(\tilde{x},\tilde{y},\tilde{z})=M\times_1W\times_2W\times_3\tilde{z},$$
for $M=\left[
  \begin{array}{cc}
    0 & \frac{1}{2} \\
    -\frac{1}{2} & 0 \\
  \end{array}
\right]$.
It follows that
$$|\langle\calC,\calC_2(\tilde{x},\tilde{y},\tilde{z})\rangle|=|\langle C_c(W,\tilde{z}),M\rangle|=\frac{\sqrt{2}}{2}\|C_c(W,\tilde{z})\|_F$$
and we conclude that maximization of $|\langle\calC,\calC_2(\tilde{x},\tilde{y},\tilde{z})\rangle|$ is equivalent to maximization of $\|C_c(W,\tilde{z})\|_F$.

Maximization of $\|C_c(W,\tilde{z})\|_F$ corresponds to multilinear rank-$(2,2,m)$ structure-preserving approximation of $\calC$. Similarly as it was done in~\cite{BeKre17} for the best antisymmetric multilinear rank-$r$ approximation, we can establish an equivalence between the best partially antisymmetric multilinear rank-$(2,2,m)$ approximation and the best unstructured rank-$1$ approximation of a partially antisymmetric tensor.

\begin{Proposition}\label{prop:phopm}
Let $\calC\in\R^{n\times n\times m}$ be a partially antisymmetric tensor. Then
\begin{align}
& \max\left\{\|\calC\times_1U^T\times_2U^T\times_3z^T\| \ : \ U\in\R^{n\times2}, U^TU=I_2, \|z\|_2=1\right\} \nonumber \\
& = \sqrt{2}\max\left\{|\calC\times_1u_1^T\times_2u_2^T\times_3z^T| \ : \ \|u_1\|_2=\|u_2\|_2=\|z\|_2=1, [u_1u_2]^T[u_1u_2]=I_2\right\} \label{phomp1} \\
& = \sqrt{2}\max\left\{|\calC\times_1v_1^T\times_2v_2^T\times_3z^T| \ : \ \|v_1\|_2=\|v_2\|_2=\|z\|_2=1\right\}. \label{phomp2}
\end{align}
\end{Proposition}

\begin{proof}
Take $\alpha=\calC\times_1u_1^T\times_2u_2^T\times_3z^T$. From the relations~\eqref{Cc} and~\eqref{Ccvalue} we see that, for every partially antisymmetric tensor $\calC$ and for $U=[u_1 u_2]$,
$$\|\calC\times_1U^T\times_2U^T\times_3z^T\|^2=\left\|\left[
  \begin{array}{cc}
    0 & \alpha \\
    -\alpha & 0 \\
  \end{array}
\right]\right\|_F^2=2\alpha^2,$$
which proves~\eqref{phomp1}.

Obviously, expression~\eqref{phomp1} is less or equal than~\eqref{phomp2}. Take the vectors $v_1$, $v_2$, $z$ that maximize~\eqref{phomp2}. There is an upper-triangular $2\times2$ matrix $R$ such that $|r_{11}|\leq1$, $|r_{22}|\leq1$ and
$$\left[
  \begin{array}{cc}
    v_1 & v_2 \\
  \end{array}
\right]=\left[
  \begin{array}{cc}
    u_1 & u_2 \\
  \end{array}
\right]R$$
is thin QR decomposition of $[v_1 v_2]$.
Using the antisymmetry in two modes we have
\begin{align*}
|\calC\times_1v_1^T\times_2v_2^T\times_3z^T| & = |\calC\times_1r_{11}u_1^T\times_2(r_{12}u_1^T+r_{22}u_2^T)\times_3z^T| \\
& = |\calC\times_1r_{11}u_1^T\times_2r_{22}u_2^T\times_3z^T| \\
& = |r_{11}r_{22}||\calC\times_1u_1^T\times_2u_2^T\times_3z^T|\leq|\calC\times_1u_1^T\times_2u_2^T\times_3z^T|
\end{align*}
This proves that the value of~\eqref{phomp1} is equal to the value of~\eqref{phomp2}.
\end{proof}

Therefore, following the previous discussion and the result of the Proposition~\ref{prop:phopm} we obtained the equivalence between our Algorithm~\ref{agm:pantisymcp} and the best unstructured rank-$1$ approximation. Then, same as in the fully antisymmetric case, the convergence result for HOPM from~\cite{Usch15} holds.

\section{Numerical experiments}\label{sec:num}

We provide numerical examples for the comparison of the CP rank-$1$ approximation with a posteriori antisymmetrization (Algorithm~\ref{agm:cp+antisym}) and the antisymmetry preserving CP (Algorithm~\ref{agm:antisymcp}). Additionally, for the sake of completeness, we compare these algorithms with the CP-ALS algorithm (Algorithm~\ref{agm:cpals}) with $r=6$, the algorithm that does not preserve antisymmetry. As we will show, antisymmetry preserving CP outperforms CP with a posteriori antisymmetrization in terms of accuracy, which was expected, but also in execution time, while CP-ALS showed to be much slower then the other two algorithms, and it also completely destroys the antisymmetric property.

All algorithms are implemented and tested in Julia programming language~\cite{bezanson2017julia}, version 1.8.1, on a personal computer, with \texttt{BenchmarkTools}~\cite{BenchmarkTools.jl-2016} package, used for determining the execution times of the algorithms (function \texttt{@btime}) and \texttt{TensorToolbox}~\cite{TensorToolbox.jl} package for tensor calculations.

For a given tensor $\calA$ and an approximation $\widetilde{\calA}$, we are looking at the relative error $\|\calA - \widetilde{\calA}\|/\|\calA\|$. We run CP-ALS algorithm, both on its own and within CP with a posteriori antisymmetrization with tolerance $10^{-8}$, and we stop antisymmetry preserving CP algorithm when either the relative error or the difference between relative errors in two consecutive iterations falls bellow $10^{-8}$.

\subsection{Example 1}\label{sec:example1}
First we generate an antisymmetric tensor $\calA$ of size $n\times n\times n$ and rank $6$, by randomly selecting three vectors $x,y,z$ of size $n$ and defining $\calA = 6\calA_6(x,y,z)$, where $\calA_6(x,y,z)$ is defined in~\eqref{antiform}. In this example we know that $\calA$ has the proposed structure. We evaluate and compare the accuracy and the speed of our algorithms. The results for different $n$ are presented in Table~\ref{tbl:rank6-ten}. The best result in each column is bolded.

\begin{table}[h!]{\small
\begin{tabular}{|r||c|c|c|c|c|c|}
	\hline
	$n$ & \multicolumn{2}{|c}{$10$} & \multicolumn{2}{|c}{$25$} & \multicolumn{2}{|c|}{$50$} \\
	\hline
	  & err & time & err & time &err & time \\
	  \hline
   CP+antisym & $0.8333$ & $224\, \mu s$  & $0.8333$ & $905.9\, \mu s$ & $0.8333$ & $\mathbf{3.983\, ms}$\\
   \hline
   antisymCP & $\mathbf{8.21\times 10^{-16}}$ & $\mathbf{69.9\, \mu s}$  & $\mathbf{1.34 \times 10^{-15}}$ & $\mathbf{502.5\, \mu s}$ & $\mathbf{1.66\times 10^{-15}}$  & $8.283\, ms$\\
   \hline
   CP-ALS  & $5.27\times 10^{-6}$ & $8.472\, ms$ & $1.998 \times 10^{-7}$ & $26.282\, ms$ & $8.43\times 10^{8}$  &  $187.625\, ms$ \\
   \hline
\end{tabular}}
\label{tbl:rank6-ten}
\caption{Evaluation of CP algorithm with a posteriori antisymmetrization (CP+antisym - Algorithm~\ref{agm:cp+antisym}), antisymmetry preserving CP (antisymCP - Algorithm~\ref{agm:antisymcp}) and CP-ALS with $r=6$ (Algorithm~\ref{agm:cpals}) in terms of relative error  $\|\calA - \widetilde{\calA}\|/\|\calA\|$ and execution times obtained by function \texttt{@btime}. }
\end{table}
	
Even though the execution time of CP with a posteriori antisymmetrization is comparable to antisymmetry preserving CP, by first approximating with non-antisymmetric tensor of rank 1, CP with a posteriori antisymmetrization loses the underlying structure, and results in an approximation with large error. CP-ALS manages to find a good non-antisymmetric approximation, but it requires much more time, so disregarding the antisymmetric property did not help either with accuracy or with execution times. Overall, antisymmetry preserving CP achieves best results.

Here, same as in the following examples, the initial vectors in Algorithm~\ref{agm:antisymcp} are taken as random vectors. If we initialize the algorithm using the HOSVD, the number of iterations decreases, but the execution time increases because of the additional time needed to perform the HOSVD.

\subsection{Example 2}
Now we construct an antisymmetric tensor by element-wise as
\begin{align*}
\mathcal{A}(i,j,k)=&\sin(x_i)\sin(y_j)\sin(z_k) + \sin(y_i)\sin(z_j)\sin(x_k)+ \sin(z_i)\sin(x_j)\sin(y_k)- \\
&\sin(y_i)\sin(x_j)\sin(z_k)- \sin(x_i)\sin(z_j)\sin(y_k)+ \sin(z_i)\sin(y_j)\sin(x_k),
\end{align*}
where $x_i$, $y_j$ and $z_k$ are sets of $n$ equidistant points on intervals $[0,1]$, $[2,10]$ and $[1,3]$, respectively. This type of tensor appear in signal processing applications. Accuracy and speed of our algorithms for different $n$ are presented in Table~\ref{tbl:antisym-func-ten}. The best result in each column is bolded.

\begin{table}[h!]{\small
\begin{tabular}{|r||c|c|c|c|c|c|}
	\hline
	$n$ & \multicolumn{2}{|c}{$10$} & \multicolumn{2}{|c}{$25$} & \multicolumn{2}{|c|}{$50$} \\
	\hline
	  & err & time & err & time &err & time \\
	  \hline
   CP+antisym & $0.8333$ & $220.9\, \mu s$  & $0.8333$ & $912.7\, \mu s$ & $0.8333$ & $\mathbf{3.95\, ms}$\\
   \hline
   antisymCP & $\mathbf{7.55\times 10^{-16}}$ & $\mathbf{111.5\, \mu s}$  & $\mathbf{9.2 \times 10^{-16}}$ & $\mathbf{3.439\, \mu s}$ & $\mathbf{1.575\times 10^{-15}}$  & $26.898\, ms$\\
   \hline
   CP-ALS  & $4.02\times 10^{-7}$ & $7.659\, ms$ & $3.91 \times 10^{-9}$ & $29.453\, ms$ & $8.492\times 10^{-7}$  &  $86.045\, ms$ \\
   \hline
\end{tabular}}
\label{tbl:antisym-func-ten}
\caption{Evaluation of CP algorithm with a posteriori antisymmetrization (CP+antisym - Algorithm~\ref{agm:cp+antisym}), antisymmetry preserving CP (antisymCP - Algorithm~\ref{agm:antisymcp}) and CP-ALS with $r=6$ (Algorithm~\ref{agm:cpals}) in terms of relative error  $\|\calA - \widetilde{\calA}\|/\|\calA\|$ and execution times obtained by function \texttt{@btime}. }
\end{table}

Similarly as in Example~\ref{sec:example1}, antisymmetry preserving CP beats two other methods in terms of accuracy and speed of getting an accurate solution.

In the next two examples we use tensors of smaller size, because ranks of those tensors increase with the size, and, since we are approximating by a rank-$6$ tensor, we want to use tensors for which it makes sense to do this type of approximation.

\subsection{Example 3}\label{sec:example3}
Now we generate an antisymmetric tensor that does not necessarily have the structure~\eqref{antiform}, by discretizing the function $f(x,y,z) = \exp(x^2+2y^2+3z^2)$ on a grid $\xi_i = (i-1)/(n-1)$, $i=1,\ldots,n$, and then applying the antisymmetrizer~\eqref{anti3d}. We test for the different values of $n$ and show the results in Table~\ref{tbl:func-ten}. Again, the best result in each column is bolded.

\begin{table}[h!]{\small
	\begin{tabular}{|r||c|c|c|c|c|c|}
		\hline
		$n$ & \multicolumn{2}{|c}{$3$} & \multicolumn{2}{|c}{$5$} & \multicolumn{2}{|c|}{$7$} \\
		\hline
		& err & time & err & time &err & time \\
		\hline
		CP+antisym & $0.8333$ & $185\, \mu s$  & $0.8339$ & $260.2\, \mu s$ & $0.8345$ & $276.6\, \mu s$\\
		\hline
		antisymCP & $\mathbf{1.61\times 10^{-14}}$ & $\mathbf{17.1\, \mu s}$  & $\mathbf{0.0557}$ & $\mathbf{87.4\, \mu s}$ & $\mathbf{0.0802}$ & $\mathbf{130.4\, \mu s}$\\
		\hline
		CP-ALS  & $2.55\times 10^{-5}$ & $9.875\, ms$ & $\mathbf{0.0557}$ & $26.282\, ms$ & $\mathbf{0.0802}$  &  $4.538\, ms$ \\
		\hline
	\end{tabular}}
	\label{tbl:func-ten}
	\caption{Evaluation of CP algorithm with a posteriori antisymmetrization (CP+antisym - Algorithm~\ref{agm:cp+antisym}), antisymmetry preserving CP (antisymCP - Algorithm~\ref{agm:antisymcp}) and CP-ALS with $r=6$ (Algorithm~\ref{agm:cpals}) in terms of relative error  $\|\calA - \widetilde{\calA}\|/\|\calA\|$ and execution times obtained by function \texttt{@btime}. }
\end{table}

Antisymmetry preserving CP achieves the best execution times. When the tensor can be well approximated by the CP approximation of the form~\eqref{antiform}, it also achieves the best accuracy ($n=3$). Otherwise, it results in the same error as CP-ALS, but much lower execution times ($n=5,7$).

\subsection{Example 4}
We generate a random tensor of size $n\times n\times n$ and antisymmetrize it with~\eqref{anti3d}. We compare the three algorithms and present the results in Table~\ref{tbl:rand-ten}. The best result in each column is bolded.

\begin{table}[h!]{\small
	\begin{tabular}{|r||c|c|c|c|c|c|}
		\hline
		$n$ & \multicolumn{2}{|c}{$3$} & \multicolumn{2}{|c}{$5$} & \multicolumn{2}{|c|}{$7$} \\
		\hline
		& err & time & err & time &err & time \\
		\hline
		CP+antisym & $0.8333$ & $184.6\, \mu s$  & $0.8546$ & $ 336.5\, \mu s$ & $0.9242$ & $743.3\, \mu s$\\
		\hline
		antisymCP & $\mathbf{3.616\times 10^{-16}}$ & $\mathbf{17\, \mu s}$  & $0.3432$ & $\mathbf{139.9\, \mu s}$ & $0.723$ & $\mathbf{493.2\, \mu s}$\\
		\hline
		CP-ALS  & $8.11\times 10^{-8}$ & $14.364 \, ms$ & $\mathbf{0.2716}$ & $20.172\, ms$ & $\mathbf{0.6393}$  &  $421.051\, ms$ \\
		\hline
	\end{tabular}}
	\label{tbl:rand-ten}
	\caption{Evaluation of CP algorithm with a posteriori antisymmetrization (CP+antisym - Algorithm~\ref{agm:cp+antisym}), antisymmetry preserving CP (antisymCP - Algorithm~\ref{agm:antisymcp}) and CP-ALS with $r=6$ (Algorithm~\ref{agm:cpals}) in terms of relative error  $\|\calA - \widetilde{\calA}\|/\|\calA\|$ and execution times obtained by function \texttt{@btime}. }
\end{table}

Similarly as in the previous example, when a tensor can be well approximated by CP decomposition with six summands (here for $n=3$), antisymmetry preserving CP achieves the best results. For $n=5,7$, antisymmetry preserving CP gives somewhat worse results than CP-ALS in terms of accuracy, but still gives the approximation in much sorter times, and CP-ALS does not preserve the antisymmetry. Note that the rank of a random antisymmetric tensor is much higher than six. This is the reason why all approximations produce high relative error.

\subsection{Example 5. Partial antisymmetry}
For the partial antisymmetry, we compare Algorithm~\ref{agm:cp+pantisym}, CP with partial a posteriori antisymmetrization, and Algorithm~\ref{agm:pantisymcp}, CP preserving partial antisymmetry, with standard CP-ALS (Algorithm~\ref{agm:cpals}) with $r=2$, which ignores the structure.

Here, regardless of how we construct the tensor $\calA$, all methods give approximately the same error. Again, the CP preserving partial antisymmetry stands out in terms of execution times. We present results in Table~\ref{tbl:partial}, with tensors $\calA_1$, $\calA_2$ and $\calA_3$ defined as follows:
\begin{itemize}
	\item $\calA_1$ is an $8\times 8\times 10$ tensor constructed by randomly selecting vectors $x, y, z$ of sizes $8,8,10$, respectively, and setting $\calA = 2\calC_2$, where $\calC_2$ is defined in~\eqref{antiform_partial}.
	\item $\calA_2$ is a $5\times 5\times 7$ tensor constructed from the function the same way as in the Example~\ref{sec:example3}.
	\item $\calA_3$ is a $5\times 5\times 4$ tensor generated by partially antisymmetrizing a tensor with randomly selected elements, using $\text{anti}_{1,2}$ operator.
\end{itemize}

\begin{table}[h!]{\small
	\begin{tabular}{|r||c|c|c|c|c|c|}
		\hline
		tensor & \multicolumn{2}{|c}{$\calA_1$} & \multicolumn{2}{|c}{$\calA_2$} & \multicolumn{2}{|c|}{$\calA_3$} \\
		\hline
		& err & time & err & time &err & time \\
		\hline
		CP+pantisym & $1.88\times 10^{-16}$ & $ 202.7\, \mu s$  & $0.1001$ & $269.4\, \mu s$ & $0.7175$ & $569.2\, \mu s$\\
		\hline
		pantisymCP & $5.832\times 10^{-16}$ & $\mathbf{30.2\, \mu s}$  & $0.1001$ & $\mathbf{53.70\, \mu s}$ & $0.7175$ & $\mathbf{106.6\, \mu s}$\\
		\hline
		CP-ALS  & $6.774\times 10^{-16}$ & $1.051 \, ms$ & $0.1001$ & $1.282 \, ms$ & $0.7175$  &  $3.026\, ms$ \\
		\hline
	\end{tabular}}
	\label{tbl:partial}
	\caption{Evaluation of CP algorithm with partial a posteriori antisymmetrization (CP+pantisym - Algorithm~\ref{agm:pantisymcp}), antisymmetry preserving partial CP (pantisymCP - Algorithm~\ref{agm:pantisymcp}) and CP-ALS with $r=2$ (Algorithm~\ref{agm:cpals}) in terms of relative error  $\|\calA - \widetilde{\calA}\|/\|\calA\|$ and execution times obtained by function \texttt{@btime}. }
\end{table}

\section{Conclusion}\label{sec:conclusion}

We described an antisymmetric tensor format $\calA_6(x,y,z)$ determined by only three vectors, $x,y,z\in\R^n$. For any $n$, tensors of the form $\calA_6(x,y,z)$ have rank at most six. We developed an ALS algorithm for structure-preserving low-rank approximation of an antisymmetric tensor $\calA$ by a tensor of the form $\widetilde{A}=\calA_6(x,y,z)$. In order to obtain our algorithm, we wrote the objective function as three different quadratic forms, given explicitly, one for each mode. The algorithm works in a way that, in each microiteration, a quadratic optimization problem for the corresponding tensor mode is solved.

We showed that our minimization problem
$$\|\calA-\calA_6(x,y,z)\|\rightarrow\min$$
for $x,y,z\in\R^n$ can be viewed as a minimization problem for orthonormal vectors $\tilde{x},\tilde{y},\tilde{z}\in\R^n$,
$$\|\calA-\calA_6(\tilde{x},\tilde{y},\tilde{z})\|\rightarrow\min.$$
Further on, we demonstrated that this minimization problem is equivalent to maximization problem
$$\|\calA\times_1V^T\times_2V^T\times_3V^T\|\rightarrow\max,$$
where $V\in\R^{n\times3}$ is a matrix with orthonormal columns. The prior maximization problem corresponds to the problem of the best multilinear low-rank approximation of antisymmetric tensors. Since antisymmetric multilinear low-rank approximation is equivalent to the best unstructured rank-$1$ approximation, we were able relate our algorithm to HOPM. Therefore, the global convergence results for HOPM apply here.

For the tensors with partial antisymmetry we established a partially antisymmetric tesnor format $\calC_2(x,y,z)$ determined by three vectors, $x,y\in\R^n$, $z\in\R^m$. Tensors of the form $\calC_2(x,y,z)$ have rank two. We created a similar ALS algorithm for structure-preserving rank-$2$ approximation of a partially atisymmetric tensor $\calC$ by a tensor of the form $\widetilde{C}=\calC_2(x,y,z)$. Analogously to the fully antisymmetric case, we verified that the algorithm in question is equivalent to HOPM.

The method described in the paper can be generalized to solve the approximation problem for different antisymmetric structures. Given that the target format can be written as a sum of multilinear terms, the underlying linearity in each mode would lead to quadratic optimization problems which would be handled in the same way, with different coefficient matrices and vectors. For example, instead of antisymmetric rank-$6$ approximation, this way, one could find antisymmetric rank-$6r$ approximation represented by $3r$ vectors. The paper limited its scope to order-$3$ tensors. For antisymmetric order-$d$ tensors, analogous rank-$d!r$ approximation would be represented by $dr$ vectors.

\section*{Acknowledgments}

The authors would like to thank the anonymous referees, whose insightful remarks improved the paper, especially for the comments regarding orthogonality of the resulting vectors and the equivalence to HOPM.

\bibliographystyle{siam}
\bibliography{CPas}

\end{document}